\documentclass[reqno]{gokova}

 \newlength{\baseunit}               
 \newcount{\numlines}                
\newcommand{\epf}{\qed \vspace{+10pt}}
\newcommand{\zed}{\mathbb{Z}}
\newcommand{\boldH}{\mathbb{H}}
\newcommand{\com}{\mathbb{C}}
\newcommand{\proj}{\mathbb P}
\newcommand{\cm}{{\mathcal{M}}}
\newcommand{\cmbar}{{\overline{\cm}}}
\newcommand{\al}{\alpha}
\newcommand{\om}{\omega}
\newcommand{\Spec}{\operatorname{Spec}}
\newcommand{\Sym}{\operatorname{Sym}}
\newcommand{\Aut}{\operatorname{Aut}}
\newcommand{\cited}{}
\newcommand{\remind}[1]{{}}

\begin{document}
\def\E{\ifmmode{\mathbb E}\else{$\mathbb E$}\fi} 
\def\N{\ifmmode{\mathbb N}\else{$\mathbb N$}\fi} 
\def\R{\ifmmode{\mathbb R}\else{$\mathbb R$}\fi} 
\def\Q{\ifmmode{\mathbb Q}\else{$\mathbb Q$}\fi} 
\def\C{\ifmmode{\mathbb C}\else{$\mathbb C$}\fi} 
\def\H{\ifmmode{\mathbb H}\else{$\mathbb H$}\fi} 
\def\Z{\ifmmode{\mathbb Z}\else{$\mathbb Z$}\fi} 
\def\P{\ifmmode{\mathbb P}\else{$\mathbb P$}\fi} 
\def\T{\ifmmode{\mathbb T}\else{$\mathbb T$}\fi} 
\def\SS{\ifmmode{\mathbb S}\else{$\mathbb S$}\fi} 
\def\DD{\ifmmode{\mathbb D}\else{$\mathbb D$}\fi} 

\renewcommand{\a}{\alpha}
\renewcommand{\b}{\beta}
\renewcommand{\d}{\delta}
\newcommand{\D}{\Delta}
\newcommand{\e}{\varepsilon}
\newcommand{\g}{\gamma}
\newcommand{\G}{\Gamma}
\newcommand{\la}{\lambda}
\newcommand{\La}{\Lambda}
\newcommand{\n}{\nabla}
\newcommand{\var}{\varphi}
\newcommand{\s}{\sigma}
\newcommand{\Sig}{\Sigma}
\renewcommand{\t}{\tau}
\renewcommand{\th}{\theta}
\renewcommand{\O}{\Omega}
\renewcommand{\o}{\omega}
\newcommand{\z}{\zeta}

\newcommand{\ben}{\begin{enumerate}}
\newcommand{\een}{\end{enumerate}}
\newcommand{\be}{\begin{equation}}
\newcommand{\ee}{\end{equation}}
\newcommand{\bea}{\begin{eqnarray}}
\newcommand{\eea}{\end{eqnarray}}
\newcommand{\bc}{\begin{center}}
\newcommand{\ec}{\end{center}}

\newcommand{\IR}{\mbox{I \hspace{-0.2cm}R}}
\newcommand{\IN}{\mbox{I \hspace{-0.2cm}N}}

\newtheorem{thm}{Theorem}[section]
\newtheorem{cor}[thm]{Corollary}
\newtheorem{lem}[thm]{Lemma}
\newtheorem{prop}[thm]{Proposition}
\newtheorem{ax}{Axiom}
\newtheorem{conj}[thm]{Conjecture}

\theoremstyle{definition}
\newtheorem{defn}{Definition}[section]

\theoremstyle{remark}
\newtheorem{rem}{\rm\bfseries{Remark}}[section]
\newtheorem*{notation}{Notation}

\newtheorem{ques}{\rm\bfseries{Question}}[section]
\newtheorem{cons}[rem]{\rm\bfseries{Construction}}
\newtheorem{exm}[rem]{\rm\bfseries{Example}}



\title{On the tautological ring of $\cmbar_{g,n}$}
\author[T. GRABER AND R. VAKIL]{Tom Graber and Ravi Vakil
}

\thanks{The first author is partially supported by an NSF Postdoctoral
Fellowship, and the second author is partially supported
by NSF Grant DMS--9970101.}

\address{Dept. of Mathematics, Harvard University, Cambridge MA USA~02138}
\email{graber@math.harvard.edu}

\address{Dept. of Mathematics, MIT, Cambridge MA USA~02139}	
\email{vakil@math.mit.edu}

\volume{7}

\maketitle

\section{Introduction}

The purpose of this note is to prove:

\begin{thm}
\remind{theorem} $R_0(\cmbar_{g,n}) \cong \Q.$
\label{theorem}  
\end{thm}

In other words, any two top intersections in the tautological ring
are commensurate.

This provides the first genus-free evidence for an audacious conjecture of Faber
and Pandharipande that the tautological ring is Gorenstein, answering
in the affirmative a question of Hain and Looijenga; see Section
\ref{gor}.

\subsection{Background on the tautological ring}
\label{background}
\remind{background}
In this section, we briefly describe the objects under consideration
for the sake of non-experts. A more detailed informal exposition of
these well-known ideas is given in \cite{pv}.

When studying Riemann surfaces of some given genus $g$, one is
naturally led to study the moduli space $\cm_g$ of such objects.  This
space has dimension $3g-3$, and has a natural compactification due to
Deligne and Mumford, \cite{dm}, the moduli space $\cmbar_{g}$ of {\em
stable} genus $g$ curves.  More generally, one can define a moduli
space of stable $n$-pointed genus $g$ curves, denoted $\cmbar_{g,n}$,
over any given algebraically closed field (or indeed over $\Spec \zed$).
We shall work over the complex numbers.  A stable $n$-pointed genus
$g$ complex curve is a compact curve with only nodes as singularities,
with $n$ distinct labeled smooth points.  There is
a stability condition: each rational component has at least 3 special
points, and each component of genus 1 has at least 1 special point.
(A special point is a point on the normalization of the component that
is either a marked point, or a branch of a node.)  This stability
condition is equivalent to requiring that the automorphism group of
the pointed curve be finite.  If a curve is stable, then a short
combinatorial exercise shows that $2g-2+n>0$.

The open subset of $\cmbar_{g,n}$ corresponding to smooth curves is
denoted $\cm_{g,n}$.  The curves of {\em compact type} (those with compact
Jacobian, or equivalently, with a tree as dual graph) form a partial
compactification, denoted $\cm_{g,n}^c$.

Even if one is interested primarily in nonsingular curves, i.e.
$\cm_g$, one is naturally led to consider the compactification,
$\cmbar_{g}$; even if one is interested primarily in unpointed curves,
one is naturally led to consider the space $\cmbar_{g,n}$ by the behavior
of the boundary, $\cmbar_{g} \setminus \cm_{g}$.

The space $\cmbar_{g,n}$ is best considered as an orbifold or
Deligne-Mumford stack; it is nonsingular and proper (i.e. compact) of
dimension $3g-3+n$, and hence has a good intersection theory.
Statements about $\cmbar_{g,n}$ translate to universal statements
about families of curves.

The cohomology of the moduli space is very interesting.  For
example, the cohomology of $\cm_g$ is the group cohomology of the
mapping class group.  There is even further structure in the Chow
ring of $\cmbar_{g,n}$ (an algebraic version of the cohomology ring).

The tautological ring, denoted $R^*(\cmbar_{g,n})$, is a subring of
the Chow ring $A^*(\cmbar_{g,n})$.
Natural algebraic constructions typically yield Chow classes lying in
the tautological ring.  
A simple description of
the tautological ring,
due to Faber and Pandharipande
(\cite{fp} Section 0.3), is the following.

Define the cotangent line class $\psi_i \in A^*(\cmbar_{g,n})$ ($1
\leq i \leq n$) as the first Chern class of the line bundle with fiber
$T_{p_i}^*(C)$ over the moduli point $[C, p_1, \dots, p_n] \in
\cmbar_{g,n}$.  The tautological system of rings $\{ R^*
(\cmbar_{g,n}) \subset A^* ( \cmbar_{g,n} ) \}$ is defined as the set
of smallest $\Q$-subalgebras satisfying the following three
properties:
\begin{enumerate}
\item[(i)]  $R^*(\cmbar_{g,n})$ contains the cotangent 
lines $\psi_1$, \dots, $\psi_n$.
\item[(ii)]  The system is closed under pushforward via all 
maps forgetting markings:
$$
\pi_*:  A^*(\cmbar_{g,n}) \rightarrow A^*(\cmbar_{g,n-1}).
$$
\item[(iii)]  The system is closed under pushforward via all gluing maps:
$$
A^* ( \cmbar_{g, n_1 \cup \{ * \} } )  \otimes_{\Q}  A^*( \cmbar_{g_2, n_2 \cup \{ \bullet \} } )
\rightarrow A^* ( \cmbar_{g_1 + g_2, n_1 + n_2} ),
$$
$$
A^*( \cmbar_{g, n \cup \{ *, \bullet \}} ) \rightarrow A^*( \cmbar_{g+1, n} ).
$$
\end{enumerate}

The Hodge bundle $\E$ is the rank $g$ vector bundle with fiber 
$H^{0}(C,\om_C)$
over the moduli point $[C, p_1, \dots, p_n] \in \cmbar_{g,n}$.  The
$\la$-classes are defined by $\la_k = c_k(\E)$ ($0 \leq k \leq g$);
these classes (and many others) also lie in the
tautological ring.

\subsection{Is the tautological ring Gorenstein?}
\label{gor}  \remind{gor}
The study of the tautological ring was initiated in Mumford's
foundational paper \cite{m}.  Recent interest in the tautological ring
was sparked by Kontsevich's proof of Witten's conjectures,
\cite{k}, which led to Faber's algorithm for computing all top intersections in
the tautological ring, \cite{f2}.  

In the early 1990's, Faber conjectured that $R^*(\cm_g)$ is of a very
special form (appearing in print much later,
\cite{f1}); it has the
properties of the $(p,p)$-cohomology of a complex projective manifold
of dimension $g-2$.  In particular, he conjectured that it is a Gorenstein ring, i.e.
$R^{g-2}(\cm_g)=\Q$, and the intersection pairing $$ R^i(\cm_g) \times
R^{g-2-i}(\cm_g) \rightarrow R^{g-2}(\cm_g) $$ is perfect (for $0 \leq
i \leq g-2$).  (Faber also gives a very specific description of top
intersections, which determines the structure of the ring.)  Key
evidence for Faber's calculation was Looijenga's beautiful proof
\cite{l} that $R^n(\cm_g)$ has dimension 0 for $n>g-2$ and at most 1
for $n=g-2$ (generated by the class of hyperelliptic curves), and
Faber and Pandharipande's subsequent argument that $\dim
R^{g-2}(\cm_g) \geq 1$ (by explicitly showing that the class of
hyperelliptic curves is non-zero, \cite{fp} equ. (8)).

Motivated by Faber's conjecture, Hain and Looijenga asked if the
tautological ring of $\cmbar_{g,n}$ satisfies Poincare duality
(\cite{hl} Question 5.5).  The first evidence would be to check if
Theorem \ref{theorem} were true; they give an argument (\cite{hl}
Section 5.1) which Faber and Pandharipande note is incomplete
(\cite{fp} Section 0.6).

More generally, Faber and Pandharipande conjecture (or speculate) that:
\begin{enumerate}
\item $R^*(\cm_g)$ is a Gorenstein ring with socle in
codimension $g-2$, \label{mg}  \remind{mg}
\item $R^*(\cm^c_g)$ is a Gorenstein ring with socle in codimension $2g-3$,
\label{mgc} \remind{mgc}
\item $R^*(\cmbar_g)$ is a Gorenstein ring with socle in codimension $3g-3$,
\label{mbarg} \remind{mbarg}
\end{enumerate} 
where $R^*$ is the tautological ring of the appropriate space, and
$\cm_g^c$ is the moduli space of curves of compact type.  There seems
no reason (or counterexample) to prevent  extensions of
these conjectures to moduli spaces of pointed curves, and indeed such
moduli spaces naturally come up inductively from the geometry of the
boundary strata of the moduli spaces of unpointed curves.

The key first step in each case is to check that the part of the
tautological ring in the highest (expected) codimension is
one-dimensional, and this note addresses the third case (further
allowing marked points).  A proof of the first step in Cases \ref{mgc}
and \ref{mbarg} was known earlier to Faber and Pandharipande
(manuscript in preparation).  However, we feel this argument (for Case
\ref{mbarg}) is worth knowing, as it is short and conceptually simple.

\begin{rem}
\label{bigchow0}  \remind{bigchow0}  
There is no reason to expect $A_0 \cmbar_{g,n}=\Q$ in general, and in
fact there is reason to expect otherwise.  For example, by a theorem
of Srinivas (\cite{s}, based on earlier work of Roitman and Mumford),
any normal projective variety $X$ with a nonzero $q$-form ($q>1$) defined
on a nonsingular open subset missing a subset of codimension at least
two has huge $0$-dimensional Chow group.  More precisely, the degree 0
dimension 0 cycle classes cannot be parametrized by an algebraic variety; in
fact the image of the natural map from $X^n$ to the degree $n$
zero-cycles modulo rational equivalence is not surjective for any $n$.
The coarse moduli scheme $\overline{M}_{g,n}$ and the fine moduli
stack $\cmbar_{g,n}$ have the same Chow group (with $\Q$-coefficients,
\cite{vistoli} p. 614), and $\overline{M}_{g,n}$ is a normal
projective variety.  For $n>1$, the spaces $\overline{M}_{1,n}$ and
$\cmbar_{1,n}$ are isomorphic away from the closed subset
corresponding to pointed curves with nontrivial automorphisms; this
set consists of the divisor $\Delta_1$ on $\overline{M}_{1,n}$
(generically corresponding to a rational curve with $n$ marked points,
attached to a genus 1 curve with no marked points) and a set of
codimension 2.  Call the corresponding divisor on the fine moduli
stack $\delta_1$.  It is not hard to show that if $n
= 11$, and hence $n \geq 11$, $\cmbar_{1,n}$ has a nonzero 11-form vanishing on
$\delta_1$; this descends to a nonzero 11-form on the open
subset of $\overline{M}_{1,n}$ corresponding to automorphism-free
curves, which extends over (the general point of) $\Delta_1$.  Hence
$\overline{M}_{1,n}$ satisfies the hypotheses of Srinivas' theorem, and has
huge Chow group.  It does not seem unreasonable to suspect that for
fixed $g>0$ and sufficiently large $n$, $\cmbar_{g,n}$ will have a
canonical form which descends to $\overline{M}_{g,n}$ (away from a
subset of codimension two).
\end{rem}

\section{The ELSV formula in Chow} 

The essential ingredient is the remarkable formula of
Ekedahl-Lando-Shapiro-Vainshtein (\cite{elsv1} Theorem 1.1, proved in
\cite{hhv} and \cite{elsv2}):
 
\begin{thm}
\label{elsvt} \remind{elsvt}
Suppose $g$, $n$ are integers ($g \geq 0$, $n \geq 1$) such that
$2g-2+n>0$, and $\al_1$, \dots, $\al_n$ are positive integers.  
Let $H^g_\al$ is the number of degree $\sum \al_i$ genus $g$ irreducible 
branched covers of $\proj^1$ with simple branching above $r$ fixed
points, branching with monodromy type $\al$ above $\infty$, and no other branching.
Then
\begin{equation*}
\label{elsvo}
H^g_\al = \frac {r!} { \# \Aut(\al)}
 \prod_{i=1}^n \frac {{\al_i}^{\al_i}} {\al_i!} 
\int_{\cmbar_{g,n}}      \frac { 1-\la_1 + \dots \pm \la_g} {\prod (1-\al_i \psi_i)}.
\end{equation*}
\end{thm}

We note that Theorem \ref{elsvt} actually 
holds in the Chow ring.  Precisely, define
the {\em Hurwitz class} $\boldH^g_\al$ to be the zero-cycle that is the sum of
the points in $\cmbar_{g,n}$ corresponding to the source curves of the
$H^g_\al$ covers, where the $n$ marked points are the preimages of
$\infty$.  (Here, $\al$ should now be considered an $n$-tuple rather than a partition; 
the points above $\infty$ are now labeled.)
This class is well-defined as all choices of the $r$
branch points in  $\proj^1$ are rationally equivalent (as
elements of $\Sym^r \proj^1$).  

\begin{prop} 
\label{elsvchowt} \remind{elsvchowt}
In $A^*(\cmbar_{g,n})$, \remind{elsvf}
\begin{equation}
\label{elsvf}
\frac 1  {r!} 
\prod_{i=1}^n \frac  {\al_i!}  {{\al_i}^{\al_i}}
 \boldH^g_\al = 
 \left[ \frac { 1-\la_1 + \dots \pm \la_g} {\prod (1-\al_i \psi_i)}\right] _{0}
\end{equation}
\end{prop}

\noindent {\em Sketch of proof.}
Our proof of Theorem \ref{elsvt} in \cite{hhv} immediately shows
that these classes are equal after pushing forward to $\cmbar _{g}$ under
the map which forgets the marked points.  To get the equality before the
pushforward, a few minor modifications to the argument are needed.
To obtain a moduli morphism to
$\cmbar_{g,n}$, virtual localization must be used on a slightly
different space. (Warning: in \cite{hhv}, $n$ is replaced by $m$.)  We
wish to count covers with $n$ {\em marked} points over $\infty$, so we
work on $\cmbar_{g,n}(\proj^1,d)$ rather than $\cmbar_g(\proj^1,d)$;
this will change the formula by $\# \Aut(\al)$.  

In \cite{hhv}, $M^k$ corresponds to the locus where the branching over
$\infty$ is (at least) $\sum (\al_i-1)$; it should be replaced by the
locus where the branching over $\infty$ is (at least) $\sum (\al_i-1)$
{\em and} the $n$ marked points map to $\infty$.  In \cite{hhv}, $M^\al$
corresponds to an irreducible component of $M^k$ where the branching
is of type $\al$; it should be replaced by the irreducible component
of $M^k$ where the branching is of type $\al$ {\em and} where the
point $p_i$ is the branch point of order $\al_i$.  This gives an
additional multiplicity of $\prod \al_i$ in $m_\al$, which is
cancelled by the same additional multiplicity in the virtual
localization formula (coming from the presence of the marked points).

\begin{rem}  The proof of \cite{elsv2} should also be easily 
adapted to show Proposition \ref{elsvchowt}, but we have not 
checked the details. \end{rem}

\subsection{} \remind{elsvchow}
\label{elsvchow}
Notice that the right side of (\ref{elsvf}) is a symmetric polynomial
$P_{g,n}$ in the variables $\al_1$, \dots, $\al_n$ (with degrees
between $2g-3+n$ and $3g-3+n$), and the coefficients are monomials in
$\psi$-classes and up to one $\la$-class.  As the coefficients of such
a polynomial can be expressed as linear combinations of the values of
the polynomial evaluated at various $n$-tuples of positive integers,
all monomials of $\psi$-classes are $\Q$-linear combinations of
Hurwitz classes.

\section{Proof of Theorem \ref{theorem}}

The proof of the theorem is now a series of reduction statements.

\subsection{} 
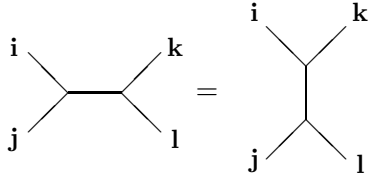
\begin{figure}
\begin{center}

\begin{picture}(150,65)(10,10)

\put(30,30){\line(-1,1){15}}
\put(30,30){\line(-1,-1){15}}
\put(30,30){\line(1,0){20}}
\put(50,30){\line(1,-1){15}}
\put(50,30){\line(1,1){15}}
 \put(7,45){\makebox(5,5){\bf i}}
 \put(7,10){\makebox(5,5){\bf j}}
 \put(68,45){\makebox(5,5){\bf k}}
 \put(68,10){\makebox(5,5){\bf l}}
 \put(80,27){\makebox(5,5){\large =}}

\put(120,40){\line(-1,1){15}}
\put(120,40){\line(1,1){15}}
\put(120,40){\line(0,-1){20}}
\put(120,20){\line(-1,-1){15}}
\put(120,20){\line(1,-1){15}}
 \put(98,58){\makebox(5,5){\bf i}}
 \put(98,2){\makebox(5,5){\bf j}}
 \put(138,58){\makebox(5,5){\bf k}}
 \put(138,2){\makebox(5,2){\bf l}}

\end{picture}

\end{center}
\caption{``Duality'' equivalence relation for trivalent
graphs \remind{soclefig}}
\label{soclefig}
\end{figure}

The smallest boundary strata (hereafter {\em top strata}) 
in $\cmbar_{g,n}$ correspond to curves whose dual graphs are trivalent
(and all components are rational).  It is well-known that they are all
rationally equivalent.  This can be seen in two ways.  First, any
point of $\cmbar _{g,n}$ corresponding to a curve with only rational
components is in the image of the standard gluing map from $\cmbar
_{0,2g+n}$ which identifies the $2g$ points in pairs.  As $\cmbar
_{0,2g+n}$ is rational, any two points on it are rationally
equivalent.  Another, purely combinatorial, method is to apply the
standard linear equivalence on $\cmbar _{0,4}$ to give explicit
rational equivalences between graphs related by the equivalence
depicted in Figure \ref{soclefig}.  A simple combinatorial exercise
shows that you can get between any two trivalent graphs using these
moves.  Thus it suffices to show that all $0$-dimensional classes in
the tautological ring are equivalent to sums of these.

\subsection{}
Top strata push forward to top strata under the natural morphisms (forgetting
markings, and gluing maps, see (ii) and (iii) in Section \ref{background}) 
between $\cmbar_{g,n}$'s.  The tautological
ring is generated by intersection of $\psi$-classes, and pushforwards
under these natural morphisms, so
it suffices to show that all top intersections of $\psi$-classes are
linear combinations of top strata.

\subsection{}  
By Section \ref{elsvchow}, as all top intersection of $\psi$-classes are
linear combinations of Hurwitz classes, it suffices to 
show that all Hurwitz classes
are linear combinations of top strata.

\subsection{}  Consider a Hurwitz cycle $\boldH^g_\al$, corresponding to 
$H^g_\al$ covers of $\proj^1$, with branching at $r+1$ given points of
$\proj^1$ (considered as an element of $\cm_{0,{r+1}}$).  Degenerate
this element to the boundary consisting of a chain of $\proj^1$'s,
with 3 special points on each component (here we use the fact that all
points of $\cmbar_{0,{r+1}}$ are rationally equivalent).  Any cover of
this configuration can have only rational components.  This is because
such a component maps to a $\proj ^{1}$ and is branched at at most
three points, and one of those branch points is a simple branch point.
It is an easy consequence of the Riemann-Hurwitz formula that such a
cover is rational.  Thus each component of the cover of the chain is
rational, with two or three special points.  The stabilization of such
a curve has only rational components with three special points, and
hence corresponds to a top stratum.
\epf

{\bf Acknowledgements:} We are grateful to R. Pandharipande and C.
Faber for introducing us to the subject, for advice and encouragement,
and for comments on an earlier version of this manuscript.  The second
author would like to thank the organizers of the Special Program on
Geometry and Topology at the Feza G\"{u}rsey Institute of Istanbul (in
May 2000), which led to many productive ongoing conversations, and
also T\"UB\.ITAK for making the Special Program possible.

\end{document}